\documentclass{amsart}

\usepackage{color} 
\usepackage{verbatim}
\usepackage{appendix}

\newtheorem{theorem}{Theorem}[section]
\newtheorem{prop}[theorem]{Proposition}

\newtheorem{cor}[theorem]{Corollary}

\newtheorem{condition}[theorem]{Condition}

\theoremstyle{definition}
\newtheorem{definition}[theorem]{Definition}

\theoremstyle{remark}
\newtheorem{remark}[theorem]{Remark}
\newtheorem{claim}[theorem]{Claim}

\newcommand{\reals}{{\bf R}}

\newcommand{\cT}{{\mathcal{T}}}

\newcommand{\cN}{{\mathcal{N}}}

\newcommand{\cR}{{\mathcal{R}}}
\newcommand{\cD}{{\mathcal{D}}}
\newcommand{\cH}{{\mathcal{H}}}

\newcommand{\cM}{{\mathcal{M}}}
\newcommand{\cP}{{\mathcal{P}}}

\renewcommand{\part}[2]{\frac{\partial #2}{\partial x_{#1}}}

\newcommand{\Ric}{{\rm{Ric}}}
\newcommand{\tr}{{\rm{tr}}}

\def\sideremark#1{\ifvmode\leavevmode\fi\vadjust{\vbox to0pt{\vss
 \hbox to 0pt{\hskip\hsize\hskip1em
 \vbox{\hsize3cm\tiny\raggedright\pretolerance10000
  \noindent #1\hfill}\hss}\vbox to8pt{\vfil}\vss}}}%
                        
\begin{document}

\title[Zero and negative eigenvalues]{Zero and negative eigenvalues of 
the conformal Laplacian}

\author[A.R. Gover]{A.R. Gover}
\address{Department of Mathematics, The University of Auckland,
Private Bag 92019, Auckland 1142, New Zealand.}
\email{r.gover@auckland.ac.nz}

\author[A. Hassannezhad]{A. Hassannezhad} 
\address{Max-Planck Institute for Mathematics
Vivatsgasse 7, 53111 Bonn, Germany}
\email{hassannezhad@mpim-bonn.mpg.de} 

\author[D. Jakobson]{D. Jakobson}
\address{Department of Mathematics and
Statistics, McGill University, 805 Sherbrooke Str. West, Montr\'eal
QC H3A 2K6, Ca\-na\-da.} \email{jakobson@math.mcgill.ca}

\author[M. Levitin]{M. Levitin}
\address{Department of Mathematics and Statistics, University of 
Reading, Whiteknights, PO Box 220, Reading RG6 6AX, United Kingdom} 
\email{M.Levitin@reading.ac.uk} 

\dedicatory{This paper is dedicated to the memory of Yuri Safarov}




\thanks{A.R.G. gratefully acknowledges support from the Royal
Society of New Zealand via Marsden Grant 13-UOA-018; 
D.J. was supported by NSERC, FQRNT and Peter Redpath fellowship.}

\begin{abstract}
We show that zero 
is not an eigenvalue of the conformal Laplacian for
generic Riemannian metrics. We also discuss non-compactness for
sequences of metrics with growing number of negative eigenvalues of
the conformal Laplacian.
\end{abstract}

\maketitle


\section{Introduction}\label{sec:intro}

In \cite{CGJP1, CGJP2}, the authors studied spectra and eigenfunctions
of conformally covariant operators on compact manifolds of dimension
$n\geq 3$.  They showed, in particular, that the number of negative
eigenvalues of conformal Laplacian is unbounded on any such manifold.
One of the questions left open in those papers was whether for a generic
Riemannian metric $g$ on a compact $n$-dimensional manifold $M$, $0$
is an eigenvalue of the conformal Laplacian $Y_g:=-\Delta_g+c_nR_g$. The operator $Y_g$ is also called the Yamabe operator.
Here $-\Delta_g$ is the nonnegative-definite Laplacian for the metric
$g$, $c_n:=(n-2)/(4(n-1))$, and $R_g$ denotes the scalar curvature of
$g$.  In this paper we address that question.

Our first main result is 
\begin{theorem}\label{thm:nozero}
For generic smooth metrics $g$ on $M$, zero is not an eigenvalue of $Y_g$.  
\end{theorem}

It follows from a transformation formula for $Y_g$ that if $0$ is an
eigenvalue of $Y_{g_0}$, then it is also an eigenvalue of $Y_{g_1}$
for all metrics $g_1$ in the conformal class $[g_0]$.  Accordingly,
one needs to change conformal class to find metrics for which $0$ is
{\it not} an eigenvalue of $Y_g$.

Also, $0$ is an eigenvalue of $Y_{g}$ for conformally {\it scalar
  flat} metrics $g$, i.e. those metrics lying in a conformal class
$[g_0]$ of a scalar-flat metric $g_0$, such that $R_{g_0}\equiv 0$.
The corresponding eigenfunction $u$ is given by $u(g_0)\equiv 1$, and
by a suitable power of the conformal factor one obtains the
eigenfunction for $g\in [g_0]$.

It is also clear that $0$ is not an eigenvalue of $Y_g$ for metrics
$g$ with {\em positive} scalar curvature $R_g$, and hence in the
corresponding conformal classes.  However, it is known that some
manifolds do not admit metrics with positive scalar curvature,
\cite{KW75}.  Accordingly, in the current paper we restrict
ourselves to metrics lying in conformal classes of metrics with {\em
  negative} scalar curvature.

\begin{remark}\label{remark:other:constants} 
We remark that our proof of Theorem \ref{thm:nozero} works 
{\em verbatim} to show an analogous
statement for operators $P_{g,c}=\Delta_g+c R_g$, where
$c\in\reals$ is any constant satisfying $c\neq 0,c\neq 1/2$.
\end{remark}
The crucial result is Proposition \ref{nonzero:derivative}, and the
only part of the proof where the numerical value of $c$ becomes
important is the argument after the equation \eqref{pde:1}, where it
is necessary that $2c-1\neq 0$, hence $c\neq 1/2$.  Also, we assume  
that the corresponding eigenfunctions are orthogonal to constants, 
hence we require that $c\neq 0$.  Note that,
$P_{g,c}$ is only conformally covariant for $c=c_n$, so for other
values of $c$ an argument using conformal perturbations (as in
\cite{BW}) should also work.  The case $c=0$ corresponds to the 
very-well studied of eigenvalues of the Laplacian $\Delta_g$.  

Perturbation theory of conformally covariant operators was previously 
considered in \cite{Canzani, Ponge}; see also \cite{Maier,Dahl03,Dahl08} 
for the corresponding results for the Dirac operator.  Applications 
to eigenvalue multiplicity of nonzero eigenvalues were discussed in 
\cite{Canzani}; arbitrary eigenvalues were considered in \cite{Ponge}, but 
the question of whether $\ker Y_g$ is generically empty was not settled.  
It seems interesting to understand whether $0$ is generically a {\em simple} 
eigenvalue of $Y_g$, among those metrics for which it {\em is} an eigenvalue 
of $Y_g$.  

Another question considered in this paper concerns the study of
sequence of metrics $g_k$ such that the number of negative eigenvalues
of $Y_{g_k}$ increases.  Recall that it was shown in \cite{CGJP1} that
the results in \cite{Lo2} imply that the number of negative
eigenvalues of $Y_g$ can become arbitrarily large for metrics $g$ on
any compact manifold of dimension $\geq 3$; thus, it seems
natural to ask what is the geometric significance of the increasing
number of negative eigenvalues of $Y_g$.

In Section \ref{sec:negative}, we show that a sequence of metrics $g_k$, such that the number 
of negative eigenvalues of $Y_{g_k}$ increases, cannot satisfy two natural ``pre-compactness'' 
conditions (see Proposition \ref{prop:non-compact}).

\section{The space of conformal structures}\label{sec:conf:structures}
Let $M$ be a compact orientable manifold of dimension $n\geq 3$; we denote by $\cM$ the 
space of all Riemannian metrics on $M$.  For simplicity, we only consider $C^\infty$ metrics 
on $M$, although the regularity can be lowered significantly.  

\begin{definition}\label{def:multiplicity:k}
Given $k\geq 1$, we denote by $\cM_{0,k}$ the 
set of all metrics $g$ on $M$ s.t. the multiplicity of $0$ as an eigenvalue of $Y_g$ is at 
least $k$.  
\end{definition}

As we remarked in section \ref{sec:intro}, if $g_0\in\cM_{0,k}$, then
so is every metric $g$ in the conformal class $[g_0]$; also that
condition is invariant under composition with diffeomorphisms of $M$.
Consider the action on $\cM$ of the group $\cP$ of (pointwise)
conformal transformations (multiplication by positive functions), as
well as by the group $\cD$ of diffeomorphisms; we shall denote by
$\cD_0$ the subgroup of $\cD$ of diffeomorphisms isotopic to
identity. It seems natural to consider the {\em Teichm\"uller space of
  conformal structures}
$$
\cT(M)\ =\ \frac{\cM/\cP}{\cD_0},
$$
or the {\em Riemannian moduli space of conformal structures}
$$
\cR(M)\ =\ \frac{\cM/\cP}{\cD},
$$
in the terminology of Fischer and Monkrief, \cite{FM96,FM97}.\footnote{If $M$
is an orientable two-dimensional manifold, then $\cT(M)$ (resp.
$\cR(M)$) are the usual Teichm\"uller (resp. moduli) spaces. In
\cite{FM97}, the space $\cT(M)$ for Haken $3$-manifolds $M$ of degree
$0$ is proposed as a configuration space for a Hamiltonian reduction
of Einstein's vacuum field equations.}  

\begin{definition}\label{def:Teich:mult:k}
We denote by $\cT_{0,k}(M)$ the Teichm\"uller space of
conformal structures corresponding to metrics $g_0\in\cM_{0,k}$,
i.e. the projection of $\cM_{0,k}$ into $\cT(M)$.
\end{definition}

The meaning of Theorem \ref{thm:nozero} is the following, and we prove this in 
Section \ref{sec:nullspace}
\begin{theorem}\label{thm:nozero:Teichmuller}
The complement $\cT_{0,1}^c$ of the set $\cT_{0,1}(M)$ in $\cT(M)$ is open and 
dense in $\cT(M)$.   
\end{theorem}


\section{Curves of metrics}\label{sec:curves}

Let $g_0$ be a metric on $M$ such that $0$ is an eigenvalue of
$Y_{g_0}$ with multiplicity $m$, so, $g_0\in\cM_{0,m}$ (recall the
definition \ref{def:multiplicity:k}). We note that it was shown in
\cite[Lemma 3.4]{BD} that the eigenvalues of $Y_g$ depend continuously
on $g$ in the $C^1$-topology (see also \cite{KS}).  Thus,
$\cM_{0,m_2}$ is a closed subset of $\cM_{0,m_1}$ for $0\leq
m_1<m_2$, in the $C^k$ topology for any $k\geq 1$. Let $C^{\omega}(I_\epsilon,\cM_{0,m})$, $\epsilon>0$ be the space of analytic curves of metrics $g(t)$, $t\in I_\epsilon=(-\epsilon,\epsilon)$ in  $\cM_{0,m}$ for $m\geq 1$.  We would like to study the space \[T_{g_0}(\cM_{0,m})=\{h\in S^2(M): \exists~ g(t)\in C^\omega(I_\epsilon,\cM_{0,m}), \epsilon>0,~ \text{with}~\dot{g}(0)=h\},\]
where $S^2(M)$ is the space of symmetric 2-tensors on $M$. 

Denote by $E_0$ the zero eigenspace of $Y_{g_0}$; it has dimension $m$.  
Let $\Pi_0$ denote the orthogonal projection into $E_0$ with respect
to $L^2(M,dV_{g_0})$.  Consider a curve $g_t$ of metrics on $M$ passing through 
$g_0$ at $t=0$; denote the $t$-derivative by $\ ^.$

Let $\dot{g}(0)=h$, i.e. $g(t)=g_0+th+o(t)$.  
We denote by $Q_{g_0,h}$ the operator 
\begin{equation}\label{projection:op}
Q_{g_0,h}:=\Pi_0\dot{Y_g}=\Pi_0(c_n\dot{R}-\dot{\Delta}):E_0\to E_0, 
\end{equation}  
Sometimes when the dependence on the metric $g_0$ is clear, we shall omit the 
subscript $g_0$ and simply write $Q_h$.

We have the following: 
\begin{prop}\label{tangent:space}
The space $T_{g_0}(\cM_{0,k})$ 
consists of all the tensors 
\begin{equation}\label{H_0:k}
\cH_{0,k}:=\{h: \ 0\ is\ an\ eigenvalue\ of\ Q_{g_0,h}\ of\ multiplicity\ \geq k \}.   
\end{equation}
\end{prop}

\noindent{\bf Proof of Proposition \ref{tangent:space}.}
We refer to \cite[page 74]{Rellich} and the discussion in \cite[\S 4,5]{Canzani} 
for basic results about the perturbation theory of conformally covariant operators; 
see also \cite{DWW}, where some important formulas that we use in our 
argument were derived.  

It follows from general theory that for a real--analytic family of
self-adjoint operators $\{Y_{g(t)}\}$, eigenvalue and eigenfunction branches 
can be chosen to depend analytically on the perturbation parameter $t$, for 
$t$ small enough, see for example \cite[Lemma 3.15]{Be73}.  Moreover, there is a positive 
constant $\epsilon$ such that, for $t$ small
enough, the number of eigenvalues of $Y_{g(t)}$ in the interval $(-\epsilon,\epsilon)$ is equal 
to $m$, where $m\ge k$ is the multiplicity of 0 for $Y_{g_0}$; and the eigenvalue derivatives
are equal to the eigenvalues of $Q_{g_0,h}$.  In particular, $Q_{g_0,h}\equiv0$ for any 
real--analytic perturbation of $g_0$ when $m=k$. 
The result is now immediate from the definition of $T_{g_0}(\cM_{0,k})$ and $\cH_{0,k}$.
  
\qed

It is well-known (see e.g. \cite{BE,FMar77}) that 
the tangent space to $\cT(M)$ at $g_0$ may be identified with the space of  all 
{\em transverse traceless} symmetric tensors $h$ satisfying  $\tr_{g_0} h=0,\delta h=0$.  
Clearly, the projection of  $\cH_{0,k}$  into the tangent space of $\cT(M)$ at $g_0$, 
consists of all transverse traceless tensors lying in $\cH_{0,k}$.  


\section{Nullspace of {$Y_g$}}\label{sec:nullspace} 
In this section we prove Theorem \ref{thm:nozero:Teichmuller}.  We keep the notation from section 
\ref{sec:curves}.  For convenience we shall assume that $g_0$ is a {\em Yamabe} metric, i.e. that 
$R_{g_0}\equiv -1$.

Theorem \ref{thm:nozero:Teichmuller} will follow from the following important result: 
\begin{prop}\label{nonzero:derivative}
There exists a {symmetric} tensor $h$ such that $Q_{g_0,h}\not\equiv 0$.  
\end{prop}
We postpone the proof that Proposition \ref{nonzero:derivative} implies Theorem 
\ref{thm:nozero:Teichmuller} until later, and first prove the Proposition.  

\noindent{\bf Proof of Proposition \ref{nonzero:derivative}.} 

Let $g_t$ be a curve of metrics real-analytic in $t$, and $\psi\not\equiv 0$ 
be an element of $E_0$ that belongs to the Rellich basis of $g_t$, 
i.e. $\psi$ is an eigenvector of the operator $Q_h$ defined in \eqref{projection:op}.

Differentiating the eigenfunction equation
$$
(-\Delta+c_n R)\psi=\lambda\psi,   
$$
we find that 
$$
\dot{\lambda}\psi=(-\Delta-\lambda+c_nR)\dot{\psi}+(c_n\dot{R}-\dot{\Delta})\psi.  
$$
It suffices to show that there exists a metric deformation $g_t$ real-analytic in $t$ such 
that $\dot{\lambda}\neq 0$.  
 
Take the inner product (with respect to $dV_{g_0}$) of both sides with $\psi$.  Since 
$(-\Delta-\lambda+c_nR)$ is self-adjoint, we find that 
$$
((-\Delta-\lambda+c_nR)\dot{\psi},\psi)=(\dot{\psi},(-\Delta-\lambda+c_nR)\psi)=0.  
$$
Then 
$$
\dot{\lambda}(\psi,\psi)=((c_n\dot{R}-\dot{\Delta})\psi,\psi).  
$$

We assume for contradiction that $\dot{\lambda}=0$ for {\em any} 
analytic perturbation $g_t$, i.e. that 
\eqref{projection:op} is identically zero for any analytic perturbation $g_t$.

We next give the expressions for $\dot{\Delta}$ and $\dot{R}$. We need to recall some notation.   
Let $C^\infty(\otimes^r T^*M)$ be the space of $(r,0)$-tensors on $M$, and $C^\infty(M)=
C^\infty(\otimes^0 T^*M)$. We consider the covariant derivative   
$$\nabla:C^\infty(\otimes^r T^*M)\to C^\infty(\otimes^{r+1} T^*M),$$
which in local coordinates is given by
$$\nabla\alpha=\sum_i\nabla_i\alpha\otimes dx_i.$$ Notice that
$d=\nabla: C^\infty(M)\to C^\infty( T^*M).$ We denote {the} formal
adjoint {of $\nabla$} by $$\delta:C^\infty(\otimes^{r+1} T^*M)\to
C^\infty(\otimes^{r} T^*M),$$ i.e. for every $\alpha\in
C^\infty(\otimes^{r} T^*M)$ and $\beta\in C^\infty(\otimes^{r+1}
T^*M),$ $(\nabla\alpha,\beta)=(\alpha,\delta\beta)$. Here,
$(\cdot,\cdot)=\int_M\langle\cdot,\cdot\rangle$, where
$\langle\cdot,\cdot\rangle$ is the pointwise inner product. We {can now} recall
  the expressions for $\dot{\Delta}$ and $\dot{R}$
computed in \cite[(2.5) and (2.6)]{DWW}, see also \cite{Be70,Be73}. 
They are,
$$
\dot{R}=-\langle h,\Ric\rangle+\delta^2h-\Delta\tr h,
$$
and 
$$
\dot{\Delta} f=-\langle h, {\nabla}^2f\rangle +\langle \delta h+\frac{1}{2} d \tr h,df\rangle, 
$$
where  $\nabla^2f$ is the Hessian of $f$, and $\delta^2$ is the formal 
adjoint of the Hessian. Recall  the pointwise inner product on $C^\infty(\otimes^{2} T^*M)$ is
$$
\langle\alpha,\beta\rangle=\sum_{i,j}\alpha^{ij}\beta_{ij};
$$ 
in particular, $\tr_gh=\langle g,h\rangle$. We know {\em apriori} that
the only metric deformations that will change the eigenvalue $0$,
{are} {\em transverse traceless} deformations $h$ of $g$.
However, we shall only insist that $h$ is {\em traceless}, and so
$\tr_g h\equiv 0$.  Then the previous expressions simplify to
$$
\dot{R}=-\langle h,\Ric\rangle+\delta^2h
$$
and 
$$
\dot{\Delta} f=-\langle h, \nabla^2 f\rangle +\langle \delta h,df\rangle.  
$$
Let $A:=\dot{\lambda}(\psi,\psi)$.  Combining the above expressions, we find that 
\begin{align*}
A&=(\langle h,\nabla^2\psi-c_n\psi\Ric \rangle -\langle \delta h,d\psi\rangle+c_n\psi\delta^2 h,\psi)\\
\; &=( h,\psi\nabla^2\psi)-c_n( h,\psi^2\Ric ) -
(  h,\nabla(\psi d\psi))+c_n(\psi\delta^2 h,\psi)\\
\; &=( h,\psi\nabla^2\psi)-c_n( h,\psi^2\Ric ) -(  h,\psi\nabla^2\psi)-
(  h,d\psi\otimes d\psi)+c_n(\psi\delta^2 h,\psi)\\
\; &=-(h,c_n\psi^2\Ric)-(  h,d\psi\otimes d\psi)+(\delta^2h,c_n\psi^2)\\
\; &=\int_M\langle h,c_n(\nabla^2\psi^2-\psi^2\Ric)-d\psi\otimes d\psi\rangle\\
\; &=\int_M\langle h,c_n\psi(2\nabla^2\psi-\psi\Ric)+(2c_n-1)d\psi\otimes d\psi\rangle.
\end{align*}

To get the last equality, we use the identity 
$\nabla^2\psi^2=2(\psi\nabla^2\psi+d\psi\otimes d\psi)$. Using the assumption that 
$\tr_g h=\langle g,h\rangle=0$, we find that 
$$
A=\int_M \langle h,c_n\psi(2\psi\mathring{\nabla}^2\psi-\psi^2\mathring{\Ric})
+(2c_n-1)(d\psi\otimes d\psi)^o\rangle,
$$
where $\mathring{V}=V-\frac{1}{n}\tr_g Vg$ is the traceless part of the 
corresponding expression $V$.

Putting $h=c_n\psi(2\psi\mathring{\nabla}^2\psi-\psi^2\mathring{\Ric})
+(2c_n-1)(d\psi\otimes d\psi)^o$ (which is symmetric and traceless), we find that $A=0$ 
if and only if 
\begin{equation}\label{pde:1}
c_n\psi(2\psi\mathring{\nabla}^2\psi-\psi^2\mathring{\Ric})
+(2c_n-1)(d\psi\otimes d\psi)^o\equiv 0.  
\end{equation}

We next remark that equation \eqref{pde:1} has no non-trivial
solutions.  Indeed, take $g_0$ to be the Yamabe metric; recall we assume that
$R_{g_0}\equiv -1$.  By assumption $\psi_{g_0}$ is an eigenfunction of
$\Delta_{g_0}$ with eigenvalue $c_n$, hence it is $L^2$-orthogonal to
the constant function, and changes sign on $M$.  
Let $\cN(\psi)$
denote the nodal set of $\psi$.  The term
$c_n\psi(2\psi\mathring{{\nabla}}^2\psi-\psi^2\mathring{\Ric})$
vanishes on $\cN(\psi)$, so it follows from \eqref{pde:1} that
$$
(d\psi\otimes d\psi)^o\equiv 0,
$$
on $\cN(\psi)$. This is equivalent to 
$$
d\psi\otimes d\psi=\frac{1}{n}|d\psi|_g^2\cdot g
$$
Now, the right-hand side has rank $n$ whenever $|d\psi|_g\neq 0$.  On the other hand, the left-hand 
side has rank $\leq 1$. The only way the equality is possible if both side are identically zero on 
$\cN(\psi)$, i.e. if 
\begin{equation}\label{sing:set}
d\psi \left|_{\cN(\psi)}\right.\equiv 0. 
\end{equation}
However, it is well-known (see e.g. \cite{Cheng,Ha,HHL,HHHN}) that the 
intersection of the nodal and critical sets of $\psi$ has locally finite Hausdorff $(n-2)$-dimensional 
measure, and so \eqref{sing:set} is impossible for non-zero $\psi$.  This 
contradiction finishes the proof 
of Proposition \ref{nonzero:derivative}. 

\qed\\

\noindent An immediate consequence of the proof of Proposition \ref{nonzero:derivative} and Proposition \ref{tangent:space} 
is the following corollary.
\begin{cor}\label{nonzeroper1}
Let $g_0$ be a Yamabe metric on $M$, $g_0\in\cM_{0,m}(M)\cap\cM_{0,m+1}(M)^c.$
Let $h=c_n\psi(2\psi\mathring{\nabla}^2\psi-\psi^2\mathring{\Ric})
+(2c_n-1)(d\psi\otimes d\psi)^o$, where $\psi=\psi_{g_0}$ is a nonzero eigenfunction of 
$Y_{g_0}$ with eigenvalue 0.   
Consider the perturbation $g(t)=g_0+th$.  Then for every $\epsilon>0$, there exists
 $|t|\le\epsilon$ such that $g(t)\in\cM_{0,m}^c$.
\end{cor}
To complete the proof of Theorem \ref{thm:nozero:Teichmuller}, we need to prove the following 
\begin{claim}\label{claim:implies}
Proposition \ref{nonzero:derivative} implies Theorem \ref{thm:nozero:Teichmuller}.  
\end{claim}
\noindent{\bf Proof of Claim \ref{claim:implies}.}  
Since $\cT_{0,1}(M)$ is a {\em closed} subspace of $\cT(M)$, its complement $\cT_{0,1}(M)^c$ 
is clearly open in $\cT(M)$. We thus need to show that $\cT_{0,1}(M)^c$ is {\em dense} in 
$\cT(M)$. We shall show that $\cM_{0,1}(M)^c$ is dense in $\cM(M)$.  
Let $g_0$ be a metric on $M$.  It suffices to show that $\cM_{0,1}(M)^c$ is dense in 
some neighbourhood $U$ of $g_0$ in $\cM(M)$.  If $g_0\in\cM_{0,1}(M)^c$, we are done, 
so we can assume that $g_0\in\cM_{0,1}(M)$.  

The proof proceeds by induction on the dimension $m$ of $E_0$.  We note that 
$m$ is finite for any $g_0$, and that Corollary \ref{nonzeroper1} was proved 
for {\em arbitrary} $m$.  
First, let $m=1$, meaning $0$ is a simple eigenvalue of $Y_{g_0}$. 
 By Corollary \ref{nonzeroper1}, we 
 know that one can choose a curve of metrics $g(t)$, real analytic in $t$, such
that $g(0)=g_0$ and  $g(t)\notin\cM_{0,1}$ for arbitrary small $t$.\footnote{It will then
  follow that (in the notation of \cite{Teytel}), the space of
  conformal structures corresponding to metrics in $\cM_{0,1}$ is of
  {\em meager codimension $1$} in the space of all conformal
  structures; we leave the details of the argument to the reader.} Hence, the proof in case
$m=1$ is complete.\\

Next, assume that we have shown that $\cM_{0,1}(M)^c$ is dense in any  
neighborhood $U$ of any metric $g_0\in \cM_{0,1}(M)$ such that 
zero is an eigenvalue of $Y_{g_0}$ with the multiplicity at most $m-1$; we would 
like to prove the corresponding statement for a metric $g_0$ such that 
$0$ is an eigenvalue of $Y_{g_0}$ with multiplicity {\em exactly} $m$.  
 By Corollary \ref{nonzeroper1}, there exists a small perturbation 
that  {\em decreases} the multiplicity $m$ of $0$ as an eigenvalue of $Y_g$.  
By the inductive hypothesis, it  follows that for any neighborhood $U$ of $g_0$, there exists a metric 
$g_1\in U,$ such that $0$ is an eigenvalue of $Y_{g_1}$ with multiplicity $\leq m-1$, 
and in a suitable neighborhood $V$ of $g_1$ (which can be chosen to satisfy 
$V\subset U$), have a nonempty intersection with $\cM_{0,1}(M)^c$. 
This completes the proof of the Claim \ref{claim:implies}, and hence also of Theorems
\ref{thm:nozero:Teichmuller} and \ref{thm:nozero}. \\

\qed



\section{Negative eigenvalues of the conformal Laplacian}\label{sec:negative} 
In \cite{CGJP1,CGJP2, El12} the authors showed that the number of negative eigenvalues of the 
conformal Laplacian cannot be uniformly bounded above on any compact 
manifold $M$ of dimension 
$n\geq 3$.  Accordingly, it seems interesting to consider sequences of metrics $g_k$ on $M$ 
where the number of negative eigenvalues of 
$Y_{g_k}=-\Delta_{g_k}+c_nR_{g_k}$ is growing.

It is known that the set of metrics $g_k$ on a manifold $M$ of dimension $n\geq 3$ is {\em pre-compact} in Gromov-Hausdorff topology if it satifies either  condition \ref{precompact:GH} or condition \ref{precompact2} below: 
\begin{condition}\label{precompact:GH}
The volume $\mathrm{Vol}(M,g_k)\leq V<\infty$ is bounded above; the injectivity 
radius $\mathrm{inj}(M,g_k)\geq r>0$ is bounded from below; the Ricci curvature $\mathrm{Ric}(M,g_k)\geq -a^2$ 
is bounded from below. 
\end{condition}
\begin{condition}\label{precompact2}
The diameter $\mathrm{diam}(M,g_k)\leq D<\infty$ is bounded above;  the Ricci curvature $\mathrm{Ric}(M,g_k)\geq -a^2$ 
is bounded from below.
\end{condition}

Consider a sequence of metrics $\tilde{g}_k$ on a fixed Riemannian  manifold such 
that the number of negative eigenvalues 
of the conformal Laplacian $Y_{\tilde{g}_k}$ goes to infinity.  It is natural to choose a 
unique Yamabe representative $g_k$ in 
the conformal class $[\tilde{g}_k]$; the scalar curvature of $g_k$ is constant and 
equal to $-1$; the number of negative eigenvalues of $Y_{g_k}$ and 
$Y_{\tilde{g}_k}$ are equal.  

\begin{prop}\label{prop:non-compact}
The sequence $g_k$ cannot satisfy the pre-compactness condition \ref{precompact:GH}; nor can 
it satisfy the condition \ref{precompact2}. 
\end{prop} 

\noindent{\bf Proof of Proposition \ref{prop:non-compact}:}
The result follows from \cite[Thm. 6.2]{Buser82} and \cite[Appendix C]{Gro}.  
Indeded, since $g_k$ is Yamabe, the number of negative eigenvalues 
of {$Y_{g_k}$} is equal to the number $N(\frac{n-2}{4(n-1)};g_k)$ of eigenvalues 
of the Laplacian $-\Delta_{g_k}$ that are less than 
$(n-2)/(4(n-1))$.  Assuming $g_k$ satisfies \ref{precompact:GH}, 
it follows from \cite[Thm. 6.2]{Buser82} that 
$N(\frac{n-2}{4(n-1)};g_k)\leq C_1<\infty$ where the constant $C_1$ only 
depends on $V,r,n,\delta$.  Similarly, assuming $g_k$ satisfies 
Condition \ref{precompact2}, it follows from Gromov's result in 
\cite[Appendix C]{Gro} that $N(\frac{n-2}{4(n-1)};g_k)\leq C_2<\infty$ 
where the constant $C_2$ only depends on $D,n,a$.  
These contradict 
the assumption on the number of negative eigenvalues of {$Y_{g_k}$}.  
\qed

Proposition \ref{prop:non-compact} shows that sequences of metrics 
with increasing number of negative eigenvalues of $Y_{g_k}$ cannot 
stay in the ``thick'' part of $\cM$ satisfying natural pre-compactness 
conditions \ref{precompact:GH} or \ref{precompact2}, and thus we cannot 
use those conditions to choose a convergent subsequence of metrics. 
On the other hand, we remark that on certain high-dimensional 
manifolds (cf. \cite{GL}) there 
exist infinitely many connected components of the set of metrics 
with {\em positive} scalar 
curvature.  Accordingly, the sequence of metrics can diverge but 
the number of negative eigenvalues 
of $Y$ can stay equal to $0$.  

 
\subsection{Example: product of a surface with another manifold.} 
We consider (a slight modification of) one of the examples discussed 
in \cite[\S 4]{CGJP1}.  Let $M$ be a manifold of dimension $d\geq 2$, and 
let $\Sigma$ be a Riemann surface of genus $\gamma\geq 2$.  Assume that $M$ admits 
a metric with positive scalar curvature, and fix 
a Yamabe metric $G$ on $M$ with scalar curvature $R_G>0$.  Fix $\epsilon>0$. 
By a result of Buser~\cite[Theorem 4]{Buser77}, for every $k\geq 1$, there 
exists a hyperbolic metric $h_k$ on $\Sigma$ such that the hyperbolic 
Laplacian $-\Delta_{h_k}$ has at least $k$ eigenvalues in the 
interval $(1/4,1/4+\epsilon)$.  Choose $k$ of those eigenvalues and denote 
them by $1/4<\lambda_{k,1}\leq\lambda_{k,2}\leq\ldots\leq
\lambda_{k,k}<1/4+\epsilon$.  Denote the corresponding eigenfunctions by 
$u_{k,j},1\leq j\leq k$.

Consider the product metric $g_k:=(G\otimes t^{-1}h_k)$ on $M\times \Sigma$, 
where $t$ is a positive constant to be chosen later.  It is easy to show that 
the scalar curvature of $g_k$ is equal to $R_G-2t$ for all $k$ (the Gauss 
curvature of $(\Sigma,h_k)$ is equal to $-1$).  If we choose 
\begin{equation}\label{t:choice1}
t>R_G/2, 
\end{equation}
then the scalar curvature of $g_k$ will be negative.  

Denote the coordinates 
on $M\times\Sigma$ by $(x,y)$.   Then the conformal Laplacian is given by 
$$
Y_{g_k}=-\Delta_{G,x}-t\Delta_{h_k,y}+\frac{d}{4(d+1)}(R_G-2t).  
$$
It follows that 
$$
Y_{g_k}u_{j,k}=\left(t\lambda_{j,k}+\frac{d(R_G-2t)}{4(d+1)}\right)u_{j,k}.  
$$
We would like to choose $t$ so that the eigenvalues 
$t\lambda_{j,k}+\frac{d(R_G-2t)}{4(d+1)}$ are all negative.  
Since $\lambda_{j,k}<1/4+\epsilon$ by assumption on $h_k$, it suffices to 
choose $t$ so that 
$$
\frac{d(2t-R_G)}{4t(d+1)}>\frac{1}{4}+\epsilon.  
$$
This can be rewritten as 
$$
\left(\frac{d}{d+1}\right)(1-R_G/2t)>1/2+2\epsilon.  
$$
A straightforward calculation shows that this is equivalent to choosing 
\begin{equation}\label{t:choice2}
t<R_G\cdot\frac{d}{d-1-4\epsilon(d+1)}.   
\end{equation}
The inequalities \eqref{t:choice1} and \eqref{t:choice2} are compatible provided 
$d/(d-1-4\epsilon(d+1))>1/2$, which is easy to achieve by choosing $\epsilon$ 
small enough.  It follows that the functions $u_{j,k}(y)$ will be eigenfunctions 
of $Y_{g_k}$ with negative eigenvalues.  

After rescaling $g_k$ by $(2t-R_G)$, we can make the scalar curvature 
$R_{g_k}\equiv -1$.  Note that the rescaling does not depend on $k$.  
It is well-known that as the number of eigenvalues of $-\Delta_{h_k}$ in 
$(1/4,1/4+\epsilon)$ increases, the injectivity radius of the metric $h_k$ 
goes to $0$, and $h_k$ leaves the ``thick'' part of the moduli space 
$\cM_\gamma$ of the hyperbolic metrics on $\Sigma.$\footnote{See e.g. \cite{Buser92}; 
for finer asymptotics of small eigenvalues we refer to \cite{Batchelor} and references 
therein.}  Accordingly, 
the injectivity radius of $(M\times\Sigma,g_k/(2t-R_G))$ also goes to $0$.  
The moduli space $\cM_\gamma$ can be compactified by adding 
surfaces with cusps; the sequence $h_k$ will then have a convergent 
subsequence in $\overline{\cM_\gamma}$, and the sequence 
$(M\times\Sigma,g_k/(2t-R_G))$ will also have a convergent subsequence.

It seems interesting to better understand under what circumstances a 
sequence of metrics $g_k$, with increasing number of negative eigenvalues of 
$Y_{g_k}$, can be made to converge in a suitable completion of the moduli 
space $\cR(M)$ of conformal structures on $M$.  

\section*{Acknowledgements}
The authors wish to thank Yaiza Canzani, Rapha\"el Ponge, Peter Sarnak and Richard 
Schoen for stimulating discussions about this problem.  The authors would like to thank the 
organizers of the workshop ``Shape optimization and spectral geometry'' at International Centre 
for Mathematical Sciences in Edinburgh, where part of this 
research was conducted. A.H. gratefully acknowledges the support and
hospitality of the Centre de recherche math\'ematiques in Montreal and the 
Max--Planck Institute for Mathematics in Bonn.  D. J. gratefully acknowledges 
the hospitality of the Department 
of Mathematics at the University College London, as well as the Department of Mathematics 
at the University of Auckland.  The authors would like to thank the anonymous referee 
for useful remarks about the paper.  


\end{document}